\documentclass[12pt,oneside,english]{amsart}
\usepackage[latin1]{inputenc}
\usepackage{amssymb}

\makeatletter
\newtheorem{theorem}{Theorem}

\newtheorem{definition}{Definition}

\newtheorem{proposition}{Proposition}

\usepackage{babel}

\makeatother
\begin{document}
\baselineskip=17pt

\title[Combinatorial identities]{Combinatorial identities generated by difference analogs of hyperbolic and trigonometric functions of order $n$}

\author{Vladimir Shevelev}
\address{Department of Mathematics \\Ben-Gurion University of the
 Negev\\Beer-Sheva 84105, Israel. e-mail:shevelev@bgu.ac.il}


\begin{abstract}
We naturally obtain some combinatorial identities finding the difference analogs of hyperbolic and trigonometric functions of order $n.$ In particular, we obtain the identities connected
with the proved in the paper the addition formulas for these analogs. \end{abstract}

\maketitle

\section{Introduction}

The original definitions of the hyperbolic and trigonometric functions of order $n$
are the following (cf. \cite{1}, point 18.2).
\begin{definition}\label{d1}
The $n$ functions
\begin{equation}\label{1}
h_s(x,n)=\frac{1}{n}\sum_{t=1}^n\omega^{(1-s)t}\exp(\omega^tx), \enskip s=1,...,n,
\end{equation}
where $\omega=\exp(\frac{2\pi i}{n}),$ are called hyperbolic functions of order $n.$
\end{definition}
In particular,
\begin{equation}\label{2}
h_1(x,1)=e^x, h_1(x,2)=\cosh x, h_2(x,2)=\sinh x.
\end{equation}
\begin{definition}\label{d2}
The $n$ functions
\begin{equation}\label{3}
k_s(x,n)=\sum_{t=0}^{\infty}\frac{(-1)^t x^{nt+s-1}}{(nt+s-1)!}, \enskip s=1,...,n,
\end{equation}
$n\geq2,$ are called trigonometric functions of order $n.$
\end{definition}
In particular,
\begin{equation}\label{4}
k_1(x,1)=e^{-x}, k_1(x,2)=\cos x, k_2(x,2)=\sin x.
\end{equation}
We consider the following equivalent definitions which could be proved directly from
Definitions \ref{d1}, \ref{d2} and the uniqueness of the solution of the Cauchy problem.
\begin{proposition}\label{p1}
a) The functions $\{h_s(x,n)\}, \enskip s=1,...,n,$ form the solution of the Cauchy problem
for the following system of ordinary differential equations
$$y_s'=y_{s-1}, \enskip s=2,3,...,n, \enskip y_1'=y_n$$
with the initials $y_1(0)=1,\enskip  y_s(0)=0,\enskip  s=2,...,n.$\newline
b) The functions $\{k_s(x,n)\},\enskip  s=1,...,n,$ form the solution of the Cauchy problem
for the following system of ordinary differential equations
$$y_s'=y_{s-1}, \enskip s=2,3,...,n, \enskip y_1'=-y_n$$
with the initials $y_1(0)=1,\enskip  y_s(0)=0, \enskip s=2,...,n.$
\end{proposition}
Note that also we have
\begin{equation}\label{5}
h_s(x,n)=\sum_{t=0}^{\infty}\frac {x^{nt+s-1}}{(nt+s-1)!}, \enskip s=1,...,n.
\end{equation}
 Proposition \ref{p1} allows to introduce the difference analogs of hyperbolic and trigonometric functions of order $n$. As usual, set $\Delta f(m)=f(m+1)-f(m).$
\begin{definition}\label{d3}
For a fixed $n$ and nonnegative integer variation $m,$ the functions $\{H_s(m,n)\}, \enskip s=1,...,n,$ are called difference hyperbolic of order $n$ if they form the solution of the following system of difference equations
\begin{equation}\label{6}
\Delta y_s(m)=y_{s-1}(m), \enskip s=2,3,...,n, \enskip \Delta y_1(m)=y_n(m)
\end{equation}
with the initials $y_1(0)=1,\enskip  y_s(0)=0, \enskip s=2,...,n.$
\end{definition}
\begin{definition}\label{d4}
For a fixed $n$ and nonnegative integer variation $m,$ the functions $\{K_s(m,n)\}, \enskip s=1,...,n,$ are called difference trigonometric of order $n$ if they form the solution of the following system of difference equations
\begin{equation}\label{7}
\Delta y_s(m)=y_{s-1}(m),\enskip  s=2,3,...,n, \enskip \Delta y_1(m)=-y_n(m)
\end{equation}
with the initials $y_1(0)=1, \enskip y_s(0)=0, \enskip s=2,...,n.$
\end{definition}
Our goal is, using the properties of functions $H_s(m,n)$ and $K_s(m,n),$ to prove
the following identities.
\begin{theorem}\label{t1}  For $m\geq0,$ we have
 \begin{equation}\label{8} 
\sum_{t\geq0}\binom {m}{nt+s-1}=\frac{1}{n}\sum_{j=1}^n
(\omega^j+1)^m\omega^{j(1-s)},\enskip s=1,...,n;
\end{equation}
 \begin{equation}\label{9}
\sum_{t\geq0}(-1)^t\binom {m}{nt+s-1}=\frac{1}{n}\sum_{j=1}^n
(\mu^{2j-1}+1)^m\mu^{(2j-1)(1-s)},\enskip s=1,...,n,
 \end{equation}
where $\mu=\exp(\frac{\pi i}{n}).$
\end{theorem}
 \newpage
 Note that formula (\ref{8}) is known (\cite{3}, \cite{4}), but formula (\ref{9}) probably is new
 (at least, it is neither in \cite{3} nor in \cite{4}). \newline
\indent Let us define the sets $\{H_s(m,n)\}, \{K_s(m,n)\}$ outside $s\in\{1,...,n\},$
putting for  $s=1,...,n,$
 \begin{equation}\label{10}
H_{-(s-1)}(m,n)=H_{n-s+1}(m,n), \enskip K_{-(s-1)}(m,n)=-K_{n-s+1}(m,n).
\end{equation}
Below we show that the definition (\ref{10}) is quite natural.
\begin{theorem}\label{t2}
(The addition formulas) For integers $m,s\geq0$ we have the identities:
\begin{equation}\label{11}
H_i(m+s, n)=\sum_{j=1}^n H_j(s,n)H_{i-j+1}(m,n), \enskip i=1,...,n;
\end{equation}
\begin{equation}\label{12}
K_i(m+s, n)=\sum_{j=1}^n K_j(s, n)K_{i-j+1}(m,n),\enskip i=1,...,n.
\end{equation}
\end{theorem}
Finally, consider circulant matrices $\mathrm{H}_n$ and $\mathrm{K}_n$
with the first row $\{(-1)^{i-1}H_i(m,n)\},\enskip i=1,...,n,$ and $\{(-1)^{i-1}K_i(m,n)\},\enskip i=1,...,n,$ respectively.
\begin{theorem}\label{t3}
$1)$ If $n$ is even, then for every $m\geq1,$  $\det{\mathrm{H}_n}=0;$\newline
$2)$ If $n$ is odd, then for every $m\geq1,$  $\det{\mathrm{K}_n}=0.$
\end{theorem}

\section{Proof of Theorem 1}

\begin{proof} Using $\Delta \binom {m}{k}=\binom {m}{k-1},$ it is easy to verify that
$H_s(m,n)$ and $K_s(m,n)$ have the following form (such that the initials evidently hold):
\begin{equation}\label{13}
H_s(m,n)=\sum_{t\geq0}\binom {m}{nt+s-1} , \enskip s=1,...,n;
\end{equation}
\begin{equation}\label{14}
K_s(m,n)=\sum_{t\geq0}(-1)^t\binom {m}{nt+s-1} , \enskip s=1,...,n.
\end{equation}
Moreover, (\ref{13}) and (\ref{14}) agree with (\ref{10}). For example, consider the equality from (\ref{10}) for $s=1$ $K_0(m,n)=-K_n(m,n).$  We have
$$K_n(m)=\binom{m}{n-1}-\binom{m}{2n-1}+\binom{m}{3n-1}-...  $$
and formally for $"s=0"$ we have
$$ K_0(m,n)=\binom{m}{-1}-\binom{m}{n-1}+\binom{m}{2n-1}-\binom{m}{3n-1}-...  .$$
Since $\binom{m}{-1}=0,$ the considered equality is evident.
 \newpage
 Furthermore, note that, by Definition \ref{d3}, \ref{d4}, $H_s(m,n),\enskip s=0,...,n$
 satisfies the difference equation $\Delta^ny-y=0,$  while $K_s(m,n),\enskip s=0,...,n$
 satisfies the difference equation $\Delta^ny+y=0.$ The characteristic equations of these
 difference equations are (cf.\cite{2})
 \begin{equation}\label{15}
\sum_{k=0}^n(-1)^{n-k}\binom {n} {k}\lambda^k\mp1 =(\lambda-1)^n\mp1=0
\end{equation}
respectively. Thus we have
 $$H_s(m,n)=\sum_{j=1}^n C_{s,j}^{(1)} (\omega^j+1)^m=$$
 \begin{equation}\label{16}
\sum_{j=1}^nC_{s,j}^{(2)}\omega^{(1-s)j} (\omega^j+1)^m, \enskip s=1,...,n.
 \end{equation}
 Further, note that to obtain $n$ distinct roots of $x^n=-1$ we could consider $\{\mu,\mu^3,...,\mu^{2n-1}\},$ where $\mu=\exp(\frac{\pi i}{n}).$  So,
 $$K_s(m,n)=\sum_{j=1}^{n}C_{s,j}^{(3)}(\mu^{2j-1}+1)^m=$$

\begin{equation}\label{17}
\sum_{j=1}^{n}C_{s,j}^{(4)}\mu^{(1-s)(2j-1)}(\mu^{2j-1}+1)^m, \enskip s=1,...,n.
 \end{equation}
Let us show that
 $$C_{s,j}^{(2)}=C_{s,j}^{(4)}=\frac{1}{n}$$
 such that
  \begin{equation}\label{18}
 H_s(m,n)= \frac{1}{n}\sum_{j=1}^n\omega^{(1-s)j} (\omega^j+1)^m, \enskip s=1,...,n;
  \end{equation}
 \begin{equation}\label{19}
K_s(m,n)=\frac{1}{n}\sum_{j=1}^{n}\mu^{(1-s)(2j-1)}(\mu^{2j-1}+1)^m, \enskip s=1,...,n.
\end{equation}
Indeed, it is easy to verify that $\Delta H_s(m,n)=H_{s-1}$ (in particular, $\Delta H_1(m,n)=H_0(m,n)=H_n(m,n));\enskip$
$\Delta K_s(m,n)=K_{s-1}$ (in particular, $\Delta K_1(m,n)=K_0(m,n)=-K_n(m,n)).$ Initials also hold, in view of identities for $s>1:$
$$\sum_{j=1}^n \omega^{j(1-s)}= 0,\enskip \sum_{j=1}^n
 \mu^{(2j-1)(1-s)}=\mu^{s-1}\sum_{j=1}^n\omega^{j(1-s)}=0.$$
 Comparing (\ref{13}) with (\ref{18}) and (\ref{14}) with (\ref{19}) we obtain (\ref{8}) and
 (\ref{9}) respectively.
\end{proof}
 \newpage  
 Using (\ref{18}), (\ref{19}) and simple transformations, we obtain, for example, formulas:
 $$K_2(m,2)=(\sqrt{2})^m \sin {\frac{\pi m}{4}},$$
 $$ H_1(m,3)=\frac{1}{3}(2^m+2\cos {\frac{\pi m}{3}}).$$
 These are A009545, A024493 \cite{5} respectively. In particular, $(\sqrt{2})^m \sin {\frac{\pi m}{4}}$ is the difference analog of $k_2(x,2)=\sin {x} .$  \newline
 Further examples: for $m\geq1,$ using (\ref{9}), we have
 $$K_1(m,5)=(2/5)(\varphi+2)^{m/2}(cos(\pi m/10) + (\varphi-1)^mcos(3\pi m/10)),$$
 $$K_2(m,5)=(2/5)(\varphi+2)^{m/2}(cos(\pi(m-2)/10) + (\varphi-1)^mcos(3\pi(m-2)/10)),$$
 $$K_3(m,5)=(2/5)(\varphi+2)^{m/2}(cos(\pi(m-4)/10) + (\varphi-1)^mcos(3\pi(m-4)/10)),$$ 
 $$K_4(m,5)=(2/5)(\varphi+2)^{m/2}(cos(\pi(m-6)/10) + (\varphi-1)^mcos(3\pi(m-6)/10)),$$
 $$K_5(m,5)=(2/5)(\varphi+2)^{m/2}(cos(\pi(m-8)/10) + (\varphi-1)^mcos(3\pi(m-8)/10)),$$
 where $\varphi$ is the golden ratio. These sequences are A289306, A289321, A289387, A289388, A289389 \cite{5}
 respectively. Note that in case $n=5$ in (\ref{9}) the third summand is 0, but if $m=0,$
 $0^0$ is accepted as 1. It is the reason why the formulas hold only for $m\geq1.$
 Note also that, using these formulas, it is easy to find all zeros of the functions $K_i(m,5).$ So, we find that  \newline
  $K_1(m,5)=0$ if and only if $m\equiv5 \pmod {10};$  \newline
  $K_2(m,5)=0$ if and only if $m=0$ or $m\equiv7 \pmod {10};$  \newline
  $K_3(m,5)=0$ if and only if $m=0, m=1$ or $m \equiv9 \pmod {10};$ \newline
  $K_4(m,5)=0$ if and only if $m=0, m=2$ or $m\equiv1 \pmod {10}; $ \newline
  $K_5(m,5)=0$ if and only if $m=0, m=1, m=2$ or $m\equiv3 \pmod {10}.$
 
 \section{Proof of Theorem 2}
 \begin{proof}
 Since the proofs for the formulas of Theorem \ref{t2} are identical, we prove
 the latter one. Using Definition \ref{d4} and (\ref{10}), let us find the values
 $K_i(1)=y_i(1)$ (Here we write $K_i(m,n)=K_i(m)$ for a fixed $n).$ Since $\Delta y_i(m)=y_{i-1}(m),$ then
\begin{equation}\label{20}
y_i(m+1)=y_i(m)+y_{i-1}(m).
 \end{equation}
 Hence, for m=0, we have $y_i(1)=0,$  except for $i=1$ and $i=2:$ $y_1(1)=1$ and
 $y_2(1)=1.$ Consequently, the sum $\sum_{j=1}^n K_j(1) K_{i-j+1}(m)$
 contains only two positive summands for $j=1,2.$  So, by (\ref{20}), we
 have
 \begin{equation}\label{21}
 K_i(m+1)=K_i(m)+K_{i-1}(m)=\sum_{j=1}^n K_j(1) K_{i-j+1}(m).
\end{equation}
\newpage
It is formula (\ref{12}) for $s=1.$ Further we use induction. Suppose, for every $m\geq0,$
 we have
 \begin{equation}\label{22}
  K_i(m+s)=\sum_{j=1}^n K_j(s) K_{i-j+1}(m), \enskip i=1,...,n.
 \end{equation}
 Then we show that
  \begin{equation}\label{23}
   K_i(m+(s+1))=\sum_{j=1}^n K_j(s+1) K_{i-j+1}(m), \enskip i=1,...,n.
  \end{equation}
  By (\ref{21}),
  \begin{equation}\label{24}
   K_i(m+s+1)=K_i((m+s)+1)=K_i(m+s)+K_{i-1}(m+s).
  \end{equation}
Further, again by (\ref{21}), the right hand side of (\ref{23}) equals
$$\sum_{j=1}^n K_j(s+1) K_{i-j+1}(m)=\sum_{j=1}^n (K_j(s)+ K_{j-1}(s))K_{i-j+1}(m)=$$
\begin{equation}\label{25}
\sum_{j=1}^n K_j(s)K_{i-j+1}(m)+\sum_{j=1}^n K_{j-1}(s)K_{i-j+1}(m)=\Sigma_1+\Sigma_2.                                                                                                             \end{equation}
According to the induction supposition (\ref{22}), we have $\Sigma_1=K_i(m+s)$ and, by                                   (\ref{24}), it is left to prove that $\Sigma_2=K_{i-1}(m+s).$ Again by the induction
supposition (\ref{22}), we have
 \begin{equation}\label{26}
  K_{i-1}(m+s)=\sum_{j=1}^n K_j(s) K_{i-j}(m).
\end{equation}
But for $\Sigma_2$ we have
\begin{equation}\label{27}
\Sigma_2=\sum_{j=1}^n K_{j-1}(s)K_{i-j+1}(m) \enskip (j-1:=j) =
\sum_{j=0}^{n-1} K_{j}(s)K_{i-j}(m).
\end{equation}
So, by (\ref{26}), (\ref{27}) and (\ref{10}) we find
$$K_{i-1}(m+s)-\Sigma_2=K_n(s)K_{i-n}(m)-K_0(s)K_i(m)=$$
$$(-K_0(s))(-K_i(m))-K_0(m)K_i(m)=0,$$
which completes the proof.
\end{proof}
For example, using (\ref{10}), for $n=3,\enskip i=1,$ we have

$$ H_1(m+s)=H_1(s)H_1(m)+H_2(s)H_3(m)+H_3(s)H_2(m), $$

$$ K_1(m+s)=K_1(s)K_1(m)-K_2(s)K_3(m)-K_3(s)K_2(m). $$

So, in particular, using (\ref{13}) and (\ref{14}), we obtain the corresponding
identities for the binomial coefficients of the form $\binom{r}{3t+i-1},\enskip r=s,m,s+m,\enskip t\geq0,\enskip i=1,2,3.$
\newpage 

\section{Proof of Theorem 3}
\begin{proof}
By the well known classic result, the determinant of a circulant matrix $\mathrm{H}$ equals
$$\prod_{t=1}^n\sum_{i=1}^n(-1)^{i-1}H_i(m,n)\omega_t^{i-1},$$
where $\{\omega_t\}, \enskip t=1,...,n,$ are all distinct roots of order $n$ from 1.
 The factor corresponding $\omega_t=1$ equals $H_1-H_2+...-H_n$ since $n$ is even. By (\ref{13}),
 we have
 $$ H_1-H_2+...-H_n=$$
 $$\binom{m}{0}+\binom{m}{n}+\binom{m}{2n}+\binom{m}{3n}+... $$
 $$-\binom{m}{1}-\binom{m}{n+1}-\binom{m}{2n+1}-\binom{m}{3n+1}-...+$$
 $$ \binom{m}{2}+\binom{m}{n+2}+\binom{m}{2n+2}+\binom{m}{3n+2}+...- $$
 $$................................................................$$
 $$-\binom{m}{n-1}-\binom{m}{2n-1}-\binom{m}{3n-1}-\binom{m}{4n-1}-... $$
 Reading over columns, we see that all consecutive  binomial coefficients occur
 with alternative signs. It is clear that $\binom {m}{m}$ occurs in the $r$-th row,
 if $m\equiv r \pmod n,\enskip 0\leq r\leq n-1,$ and other summands are zeros.
 So, for $m\geq1,$ $ H_1-H_2+...-H_n=\sum_{l=0}^m (-1)^{l}\binom{m}{l}=0$ and also $\det\mathrm{H}=0.$
 Analogously, using (\ref{14}), for odd $n$ we find that $ K_1-K_2+K_3-...+K_n=0$ and so also
 $\det\mathrm{K}=0.$
 
\end{proof}


\begin{thebibliography}{5}

\bibitem {1} Higher Trancendental Functions, Bateman Manuscript Project, Vol. 3,
ed. A. Erdelyi, 1983.
\bibitem {2} A. O. Gelfond, Calculus of Finite Differences, Nauka, Moscow, 1967 (in Russian).
\bibitem {3} Combinatorial Identities: Table $III:$ from the seven unpublished manuscripts
of H. W. Gould, ed. J. Quaintance, 2010.
\bibitem {4} M. Merca, On some power sums of sine and cosine, The Amer. Monthly
121 (2014) no.3.1, 244-248.
\bibitem {5} N.\enskip J.\enskip A.\enskip Sloane,\enskip\slshape The On-Line
Encyclopedia of Integer Sequences \upshape \enskip http://oeis.org.

\end{thebibliography}
 \end{document}